\documentclass[a4paper,12pt]{article}
\usepackage{amsmath}
\usepackage{amsfonts}
\usepackage{amsthm}
\usepackage{graphicx}

\newtheorem{lemma}{Lemma}[section]
\newtheorem{propos}[lemma]{Proposition}

\newtheorem{theorem}[lemma]{Theorem}

\newtheorem{corollary}[lemma]{Corollary}
\newtheorem{defin}[lemma]{Definition}

\newcommand{\nat}{\ensuremath{\mathbb{N}}}

\newcommand{\integers}{\ensuremath{\mathbb{Z}}}

\begin{document}


\begin{center} {\Large {\bf Meadows and the equational specification of division}}\\
    \baselineskip 13pt{\ } \vskip 0.2in
J A Bergstra\footnote{Email: j.a.bergstra@uva.nl}
\\{\ }\\ Informatics Institute, University of Amsterdam,  \\ Science Park 403, 
		1098 SJ  Amsterdam,  The Netherlands
\\{\ }\\ Y Hirshfeld\footnote{Email: joram@post.tau.ac.il}
\\{\ }\\ Department of Mathematics, Tel Aviv University, \\Tel Aviv 69978, Israel
\\{\ }\\ J V Tucker\footnote{Email: j.v.tucker@swansea.ac.uk}
\\{\ }\\ Department of Computer Science,  Swansea University, \\ Singleton Park,  Swansea, SA2 8PP, United Kingdom
\end{center}

\bigskip

\begin{abstract}
\noindent
The rational, real and complex numbers with their standard operations, including division, are partial algebras specified by the axiomatic concept of a field. Since the class of fields cannot be defined by equations, the theory of equational specifications of data types cannot use field theory in applications to number systems based upon rational, real and complex numbers. We study a new axiomatic concept for number systems with division that uses only equations: a {\it meadow} is a commutative ring with a total inverse operator satisfying two equations which imply $0^{-1} = 0$.  All fields and products of fields can be viewed as meadows. 
After reviewing alternate axioms for inverse, we start the development of a theory of meadows. We give a general representation theorem for meadows and find, as a corollary, that the conditional equational theory of meadows coincides with the conditional equational theory of zero totalized fields. We also prove representation results for meadows of finite characteristic.
\end{abstract}
\bigskip

\section{Introduction}

At the heart of the theory of data types are the ideas of specifying the properties of data using equations and conditional equations, performing calculations and reasoning using term rewriting, and modelling all data representations and implementations using algebras. The theory combines mathematical simplicity, beauty and usefulness, especially when using equations and total operations. Confidence in the scope and explanatory power of the theory was established in its first decade, 1975-85, when it was proved that {\em any} computable data type possesses a range of equational specifications with desirable properties, such as having few equations (e.g., \cite{Bergstra&Tucker.82, Bergstra&Tucker.83, Bergstra&Tucker.87}), or equations with valuable term rewriting properties (e.g., \cite{Bergstra&Tucker.95}). Since every computable data type can be equationally specified - and, indeed, there are special specifications that define all and only computable data types - we expect that any data type arising in computing can be specified by equations and studied using the theory. The search for, and study of, equational specifications of {\it particular} computational structures is long term activity, contributing to foundational thinking in diverse areas of computer science, such as programming languages, hardware verification, graphics, etc. For the theoretician, it is a challenge to develop and perfect the properties of specifications far beyond those delivered by the general theory.  

Despite achievements in many areas, one does not have far to look for a truly fundamental challenge. Algebras of rational, real and complex numbers make use of operations whose primary algebraic properties are captured by the axioms of the concept of {\it field}. The field axioms consist of the equations that define commutative rings and, in particular, two axioms that are {\it not} equations that define the inverse operator and the distinctness of the two constants. Now, division is a partial operation, because it is undefined at $0$, and the class of fields cannot be defined by any set of equations. Thus, the theory of equational specifications of data types cannot build on the theory of fields; moreover, data type theory has rarely been applied to number systems based upon rational, real and complex numbers. However, we know that, say, the field of rational numbers is a computable data type - arguably, it is the most important data type for measurement and computation. Therefore, thanks to general theory, computable data types of rational, real and complex numbers with division do have equational specifications. This fact leads to two problems: we must search for, and study, 

1. equational specifications of {\it particular} algebras of rational, reals and complex numbers with division; and, ideally,

2.  equational specifications of classes of number algebras with division that are as elegant and useful as the theory of fields. 

Having begun to tackle Problem 1 in \cite{Bergstra&Tucker.05a, Bergstra&Tucker.05b, Bergstra.06}, this paper considers Problem 2 and introduces a new axiomatisation for number systems with division, called the {\ meadow}, which uses only equations.  

A {\em meadow} is a commutative ring with unit equipped with a total unary operation $x^{-1}$, 
named inverse, that satisfies these additional equations:
\begin{eqnarray}
(x^{-1})^{-1}& = & x \\
x \cdot (x \cdot x^{-1}) & = & x. 
\end{eqnarray}
The first equation we call {\em Ref,} for {\em reflection,} and the second equation  $Ril,$ for {\em restricted inverse law}.

Meadows provide a mathematical analysis of division which is more general than the classical theory of fields. Meadows are total algebras in which, necessarily, $0^{-1} = 0$. We have used algebras with such {\it zero totalized division} in developing elementary algebraic specifications for several algebras of numbers in our previous papers \cite{Bergstra&Tucker.05a, Bergstra&Tucker.05b, Bergstra.06}. The raison d'\^{e}tre of meadows is to be a tool that extends our understanding and techniques for making specifications. Clearly, since meadows are commutative rings they also have pure mathematical interest.

Let us survey our results. In \cite{Bergstra&Tucker.05a}, an equational specification under initial algebra semantics of the zero totalized field of rational numbers was presented, and specifications for other zero totalized fields were developed in \cite{Bergstra&Tucker.05b} and \cite{Bergstra.06}. In  \cite{Bergstra&Tucker.05a} meadows were isolated by exploring alternate {\em equational} axioms for inverse. Specifically, 12 equations were found; a set $\mathit{CR}$ of 8 equations for commutative rings was extended by a set 
  $\mathit{SIP}$ of 3 equations for inverse, including {\em Ref}, and by $\mathit{Ril}$. 
  The  single sorted finite equational specification $\mathit{CR} + \mathit{SIP} + \mathit{Ril}$ has all zero totalized fields among its models and, in addition, a large class of structures featuring zero divisors.  A model of 
  $\mathit{CR} + \mathit{SIP} + \mathit{Ril}$ was baptized a {\it meadow} in 
  \cite{Bergstra&Tucker.05a}. 
Because meadows are defined by equations, finite and infinite products of zero totalized fields are meadows as well.

Our first result will be that two of the equations from $\mathit{CR} + \mathit{SIP} + \mathit{Ril}$ can be derived from the other ones. This establishes the subset $\mathit{Md,}$ consisting of 10 equations of the 12 equations,
including the 8 equations for $\mathit{CR}$ and the equations {\em Ref}  and {\em Ril} mentioned earlier. Our second result makes an intriguing  connection between meadows and commutative von Neumann regular rings.

Our main task is to start to make a classification of meadows up to isomorphism. We prove the following general representation theorem:\\

\noindent {\bf Theorem} {\em  Up to isomorphism, the non-trivial meadows are precisely the subalgebras of products of zero totalized fields.}\\

From this theorem we deduce this corollary:\\

\noindent {\bf Theorem} {\em  The equational theory of meadows and the equational theory of fields with zero totalized division are identical.}\\

This strengthens a result for closed equations in \cite{Bergstra&Tucker.05a}. Now we prove the following extension:\\

\noindent {\bf Theorem} {\em  The conditional equational theory of meadows and the conditional equational theory of fields with zero totalized division are identical.}\\
 
Next, we examine the relationship between fields and meadows of finite characteristic. The characteristic of a meadow is the smallest natural number $n \in \nat$ such that $n.1 = 1 + 1 + \ldots + 1 = 0$. A prime meadow is a meadow without a proper submeadow and without a proper non-trivial homomorphic image.

Given a positive natural number $k$, and writing $\underline{k}$ for the numeral for $k$, we can define $\mathit{Md}_{k}$ for the initial algebra of $\mathit{Md} + \{\underline{k}=0\}$, i.e., 

\begin{center}
$\mathit{Md}_{k} \cong I(\Sigma, \mathit{Md} \cup \{\underline{k}=0\}).$
\end{center}
The following results are obtained:\\

\noindent {\bf Theorem} {\em   For $k$ a prime number, $\mathit{Md}_{k}$ is the zero totalized prime field  of characteristic $k$}. \\

\noindent {\bf Theorem} {\em   For $k$ a square free number, 
$\mathit{Md}_{k}$ has cardinality $k$}. \\

In the matter of Problem 1 above, only recently, Moss  found in \cite{Moss.01} that there exists an equational specification of the ring of rationals (i.e., without division or inverse) with just {\it one} unary hidden function. In \cite{Bergstra&Tucker.05a} we proved that there exists a finite equational specification under initial algebra semantics, {\it without} hidden functions, but making use of an inverse operation, of the field of rational numbers. In  \cite{Bergstra&Tucker.05b}, the specification found for the rational numbers was extended to the complex rationals with conjugation, and in \cite{Bergstra.06} a specification was given of the algebra of rational functions with field and degree operations that are all total.
Full details concerning the background of this work can be found in  \cite{Bergstra&Tucker.05a}.

We assume the reader is familiar with the basics of ring theory 
(e.g., \cite{McCoy.64, Stoltenberg-Hansen&Tucker.99}),
algebraic specifications (e.g.,
\cite{Wirsing.90}), universal algebra (e.g., \cite{Wechler.92, Meinke&Tucker.92}) and term rewriting (e.g., \cite{Terese.03}).

\section{Axioms for fields and meadows}\label{ANA}
We will add to the axioms of a commutative ring various alternative axioms for dealing with inverse and division. The starting point is a signature $\Sigma_{CR}$ for commutative rings with unit:
\newline

{\bf signature}  $\Sigma_{CR}$

{\bf sorts}     $ring$

{\bf operations}  

$0 \colon  \to ring$;

$1 \colon  \to ring$;

$ + \colon  ring \times ring \to ring$;

$ - \colon  ring \to ring$;

$ \cdot \colon  ring \times ring\to ring$

{\bf end}\\

To the signature  $\Sigma_{CR}$ we add an inverse operator  $^{-1}$ to form the primary signature 
$\Sigma$, which we will use for both fields and meadows:\\

{\bf signature}  $\Sigma$

{\bf import}  $\Sigma_{CR}$

{\bf operations}

$ ^{-1}  \colon ring \to ring$

{\bf end}

\subsection{Commutative rings and fields}

The first set of axioms is that of a {\it commutative ring with $1$}, which establishes the standard properties of $+$, $-$, and $\cdot$. 
\newline

{\bf equations} {\it $CR$}
\begin{eqnarray}
(x + y) + z & = & x + (y + z)    \\
x + y & = & y + x  \\
x + 0 & = & x  \\
x + (-x) & = & 0  \\
(x \cdot y) \cdot z & = & x \cdot (y \cdot z)  \\
x \cdot y & = & y \cdot x    \\
x \cdot 1 & = & x          \\
x  \cdot (y + z) & = & x \cdot y + x \cdot z
\end{eqnarray}

{\bf end}
\newline

These axioms generate a wealth of properties of $+, -, \cdot$ with which we will assume the reader is familiar. We will write $x-y$ as an abbreviation of $x + (-y)$.

\subsubsection{Axioms for meadows}
Having available an axiomatization of commutative rings with unit (such as the one above), we define the equational axiomatization of meadows by 
\[\mathit{Md} = (\Sigma,CR + \mathit{Ref} + \mathit{Ril}).\]

\subsubsection{Axioms for fields}
\label{Iel}
On the basis of the axioms $CR$ for commutative rings with unit there are different ways to proceed with the introduction of division. The orthodoxy is to add the following two axioms for fields: let $Gil$ ({\em general inverse law}) and  $Sep$ ({\em separation axiom}) denote denote the following two axioms, respectively:
\begin{eqnarray}
x \neq 0 \implies x \cdot x^{-1} &=& 1 \\
0 &\neq& 1
\end{eqnarray}
Let $(\Sigma, T_{field})$ be the axiomatic specification of fields, where
$T_{field} = CR + Gil + Sep$.
About the status of $0^{-1}$ these axioms say nothing. This may mean that the inverse is: 

(1) a partial function, or 

(2) a total function with an unspecified value, or 

(3) omitted as a function symbol but employed pragmatically as a useful notation in some ``self-explanatory'' cases.

Case 3  arises in another approach to axiomatizing fields, taken in many text-books, which is
 not to have an operator symbol for the inverse at all and to add an axiom 
 \textit{Iel} (\emph{inverse existence law}) as follows:
\[x \neq 0 \implies \exists y (x \cdot y = 1).\]

Each  $\Sigma$ algebra satisfying $T_{field}$ also satisfies $\mathit{Iel}$. 
In models of $(\Sigma_{CR}, CR+\mathit{Iel}+Sep)$ the inverse is implicit as a single-valued definable relation, so we call this theory the {\it relational theory of fields} $\mathit{RTF}$.

\subsubsection{Totalized division in fields}\label{Algebra_satisfying_the_Specifications}
In field theory, if the decision has been made to use a function symbol for inverse the value of $0^{-1}$ is either left undefined, or left unspecified. However, in working with elementary specifications, which we prefer, operations are total. This line of thought leads to totalized division.

The class $Alg(\Sigma,T_{field})$ is the class of all possible \emph{total} algebras satisfying the axioms in $T_{field}$.  For emphasis, we refer to these algebras as \emph{totalized fields}.

Now, for all totalized fields $A \in Alg(\Sigma,T_{field})$ and all $x \in A$, the inverse $x^{-1}$ is defined. Let $0_A$ be the zero element in $A$. In particular,  $0_A^{-1}$ is defined.  The actual value $0_A^{-1} = a$ can be anything but it is convenient to set $0_A^{-1} = 0_A$ (see \cite{Bergstra&Tucker.05a}, and compare, e.g., Hodges  \cite{Hodges.93}, p. 695). 

\begin{defin}A field $A$ with $0_A^{-1} = 0_A$ is called {\it zero totalized}.
\end{defin} 
This choice gives us a nice equation to use,  the {\em zero inverse law} $Zil$:
\[0^{-1} = 0.\]
With $\mathit{ZTF}$, an extension of $T_{field}$, we specify the class of zero totalized fields:
\begin{center}
$\mathit{ZTF} = T_{field} + Zil = \mathit{CR}+Gil+Sep +Zil$.
\end{center}
Let $Alg(\Sigma,\mathit{ZTF})$ denote the class of all zero totalized fields.

\begin{lemma}
Each $\Sigma_{CR}$ algebra satisfying $CR+Iel+Sep$ can be expanded to a $\Sigma$ algebra with a unique inverse operator that satisfies $\mathit{ZTF}$.
\end{lemma}

\begin{proof}
To see this notice that if $x \cdot y = 1$ and $x \cdot z = 1$ it follows by subtraction of both equations that $x \cdot (y-z) = 0$.  Now:
\begin{center}
$y-z = 1 \cdot (y-z) = (x\cdot y) \cdot (y-z) = x \cdot (y-z)\cdot y = 0\cdot y = 0$,
\end{center}
which implies that $y = z$ and that the inverse is unique. Let $x^{-1}$ be the function that produces this unique value (for non-zero arguments). Choose $0^{-1}$ to be $0$ and a zero totalized field has been built. 
\end{proof}

\subsubsection{Equations for zero totalized division}
Following \cite{Bergstra&Tucker.05a}, one may replace the axioms $Gil$ and $Sep$ by other axioms for division, especially, the three equations in a unit called $SIP$ for {\it strong inverse properties}.
They are considered ``strong" because they are equations involving $^{-1}$ {\it without any guards}, such as $x\neq 0$. These three equations were used already by  Harrison in  \cite{H1998}.
\newline

{\bf equations} {\it $SIP1, SIP2\,\, \mathrm{and}\,\, SIP3$}
\begin{eqnarray}
(-x)^{-1} & = &  -(x^{-1})  \\
(x \cdot y)^{-1} & = & x^{-1} \cdot y^{-1} \\
(x^{-1})^{-1} & = & x
\end{eqnarray}

{\bf end}\\

The following was proven  in \cite{Bergstra&Tucker.05a}:
\begin{propos} 
$CR \cup SIP  \vdash  0^{-1} = 0$.
\end{propos}

\subsection{Meadows and $Ril$}\label{Meadows}

In \cite{Bergstra&Tucker.05a} we add to $CR + SIP$ the equation $Ril$ (\emph{restricted inverse law}):
\begin{center}
$x \cdot (x \cdot x^{-1})  =  x $
\end{center}
which, using commutativity and associativity, expresses that $x \cdot x^{-1}$ is $1$ in the 
presence of  $x$. 
We may write $x\cdot x^{-1}$ as $1_x$,
in which case we have the following alternative formulations of $Ril$,
 \[1_x\cdot x=x\,\, \mathrm{and} \,\, 1_x\cdot x^{-1}=x^{-1},\]
and also
$1_x=1_{x^{-1}}$.
Following \cite{Bergstra&Tucker.05a} we define:

\begin{defin} A model of $CR + SIP + Ril$ is called a {\em meadow}. 
\end{defin}

Shortly, we will demonstrate that this definition is equivalent to the definition of  a meadow given in the introduction. 
A meadow satisfying $Sep$ is called {\it non-trivial}. \\

\noindent {\bf Example} All zero totalized fields are clearly non-trivial meadows but not conversely. 
In particular, the zero totalized prime fields $\integers_p$ of prime characteristic are meadows. That the initial algebra of $CR+SIP+Ril$ is not a field follows from the fact that $(1+1) \cdot (1+1)^{-1} = 1$ cannot be derivable because it fails to hold in the prime field $\integers_2$ of characteristic 2 which is a model of these equations as well.\\

Whilst the initial algebra of $\mathit{CR}$ is the ring of integers, we found in  \cite{Bergstra&Tucker.05a} that 
\begin{lemma} The initial algebra of $\mathit{CR}+\mathit{SIP}+\mathit{Ril}$ is a computable algebra but it is not an integral domain.
\end{lemma}

\subsection{Derivable properties of meadows}
We will now derive some equational facts from the specification $\mathit{Md}$  or relevant subsets of it.

\begin{propos} ~\\ 
$\mathit{CR} + Ril \vdash x \cdot x^{-1} = 0 \leftrightarrow x = 0$. 
\end{propos}
\begin{proof}
Indeed, we have $x \cdot x^{-1} = 0 \Longrightarrow x \cdot x^{-1} \cdot x = 0 \cdot x$, by multiplication. Thus, $x=0$ by applying $Ril$ to the LHS and simplifying the RHS. The other direction is immediate from $0 \cdot x = 0$.
\end{proof}

To improve readability we denote $x^{-1}$ by $\overline{x}$ and use
$1_x= x\cdot x^{-1}$. Recall that 
$1_x=1_{\overline{x}}$.

\begin{propos}\label{XYZ} Implicit definition of inverse:\\
$\mathit{CR} + Ril  \vdash x \cdot y = 1 \rightarrow x^{-1} = y$
\end{propos}
\begin{proof} $\overline{x} = 1 \cdot \overline{x} = x \cdot y \cdot \overline{x} = 1_x \cdot y =
(1_x + 0) \cdot y = (1_x + 0 \cdot \overline{x}) \cdot y = (1_x + (x - x) \cdot \overline{x}) \cdot y =
(1_x + (x \cdot 1 - x \cdot x \cdot \overline{x}) \cdot \overline{x}) \cdot y = 
(1_x + (x \cdot x \cdot y - x \cdot x \cdot \overline{x}) \cdot \overline{x}) \cdot y =
(1_x + x \cdot x \cdot (y - \overline{x}) \cdot \overline{x}) \cdot y = 
(1_x + x \cdot (y - \overline{x})) \cdot y = 
(1_x + x \cdot y - x \cdot \overline{x}) \cdot y = x \cdot y \cdot y = 1 \cdot y = y
$
\end{proof}

\begin{propos}\label{SIP12} Derivability of SIP1 and SIP2:\\
1.   $\mathit{Md}  \vdash  (xy)^{-1}=x^{-1}y^{-1}$\\
2.   $\mathit{Md}  \vdash  (-x)^{-1}=-(x^{-1})$
\end{propos}

\begin{proof}

1. First we show that $1_{x y}=1_x\cdot 1_y$.
Indeed we have:
 $1_{x y}\cdot 1_x\cdot 1_y= x\cdot y\cdot \overline{xy}\cdot
x\cdot  \overline{x} \cdot y\cdot \overline{y}$
Applying $Ril$ twice we have $x\cdot y\cdot x\cdot \overline{x}
\cdot y\cdot \overline{y}=x\cdot y$ , and therefore
$1_{x y}\cdot 1_x\cdot 1_y= x\cdot y\cdot \overline{xy} =1_{x y}$.
On the other hand applying $Ril$ once we have $x\cdot y\cdot
\overline{xy}\cdot x\cdot y= x\cdot y $ and therefore
$1_{x y}\cdot 1_x\cdot 1_y= x\cdot y\cdot \overline{x}\cdot
\overline{y}=1_x\cdot 1_y$
This proves the auxiliary equation. Now:
$\overline{xy}=\overline{xy}\cdot 1_{xy}=\overline{xy}\cdot
1_x\cdot 1_y=\overline{xy}\cdot x\cdot\overline{x}\cdot
  y\cdot\overline{y}=1_{xy}\cdot\overline{x}\cdot\overline{y}=1_{x}\cdot
  1_y\cdot\overline{x}\cdot\overline{y}=\overline{x}\cdot\overline{y}.$

2. The fact that $\overline{-1}=-1$ follows by an application of Proposition \ref{XYZ}
to  $(-1)\cdot(-1)=1$ which is a consequence of $\mathit{CR}$.
We now conclude  with the help of 1:
$\overline{-x}=\overline{(-1)\cdot x}=\overline{(-1)}\cdot
\overline{x}=(-1)\cdot \overline{x}=- \overline{x}$
\end{proof}

Thanks to Proposition \ref{SIP12} we obtain:
\begin{corollary}
$\mathit{Md}$ axiomatizes the meadows, i.e. $\mathit{Md}$ is equivalent to
$\mathit{CR}+\mathit{SIP}+\mathit{Ril}$.
\end{corollary}

\begin{propos}  ~\\
1.   $\mathit{CR} + Ril + \mathit{SIP2} \vdash x^2 = x \rightarrow x = x^{-1}$\\
2.   $\mathit{Md}\vdash x^3 = x \rightarrow x = x^{-1}$, and \\
3.   $\mathit{Md}\vdash x^4 = x \rightarrow x = x^{-2}$.
\end{propos}
\begin{proof}~\\
1. $x= x \cdot x \cdot x^{-1} = x \cdot x^{-1}= x \cdot (x \cdot x)^{-1} = 
x \cdot x^{-1} \cdot x^{-1} = x^{-1}.$~\\
2. From the assumption we obtain $x^3 \cdot x^{-1} = x \cdot x^{-1}$ and then
$x \cdot x = x \cdot x^{-1}$. Thus $x \cdot x \cdot x^{-1} = x \cdot x^{-1} \cdot x^{-1} $ whence
 $  x = ((x \cdot x^{-1} \cdot x^{-1})^{-1})^{-1} = (x^{-1} \cdot x \cdot x)^{-1} = x^{-1}.$~\\
3. From the assumption we obtain $x^4 \cdot x^{-1} = x \cdot x^{-1}$ and then $x^3  = x \cdot x^{-1}$,
from which we get $x^3 \cdot x^{-1} = x \cdot x^{-1}\cdot x^{-1}$ and $x^2 = x^{-1}$.
\end{proof}

\subsection{Meadows and von Neumann regular rings with unit}
A {\em commutative von Neumann regular ring} (e.g., see \cite{McCoy.64, Good.79}) is a $\Sigma_{CR}$ algebra that satisfies $\mathit{CR}$ and which in addition satisfies the following axiom
{\em regular ring} ({\em RR}):
\[\forall x .\exists y. (x \cdot y \cdot x  =   x).  \]
A value $y$ which satisfies $x \cdot y \cdot x  =   x$ is called a {\em pseudoinverse} of $x$. 

Because $\mathit{Ril}$ indicates that $x^{-1}$ is a pseudoinverse of $x$, the $\Sigma_{CR}$-reduct of a meadow is a commutative von Neumann regular ring and {\em  every meadow is an expansion of a von Neumann regular ring}. As it turns out a converse is true. We acknowlege Robin Chapman (Exeter UK) for pointing out to us the following observation:
\begin{lemma}\label{vNrrexpansion}
Every commutative regular von Neumann ring can be expanded to a meadow. Moreover, this expansion is unique.
\end{lemma}
First, we notice a lemma that holds for any commutative ring.
\begin{lemma}\label{uniqueinv}Given an $x$, any $y$ with $x \cdot x \cdot y =  x$ and $y \cdot y \cdot x = y$ is unique. 
\end{lemma}
\begin{proof}Assume that, in addition, $x \cdot x \cdot z  =   x$ and $z \cdot z \cdot x = z$. By subtracting the first equations of both pairs, we get $x \cdot x \cdot(y-z)  =  0$, which implies $x \cdot x \cdot(y-z)\cdot y  =  0 \cdot y$, on multiplying both sides by $y$. Since  $x \cdot x \cdot y =  x$, we deduce that $x \cdot(y-z)  =  0 $ and that $x \cdot y  =  x \cdot z$. Now, substituting into $y \cdot y \cdot x = y$, this yields $ y \cdot z \cdot x = y$; and substituting into $z \cdot z \cdot x = z$ it yields $z \cdot y \cdot x = z$; taken together, we conclude $y = z$. 
\end{proof}
\begin{proof}
Then we proceed with the proof of Lemma \ref{vNrrexpansion}.
Suppose that $\Sigma_{CR}$ algebra $A$ satsifies $\mathit{RR}$. First, expand the $A$ to an algebra $A^{\prime}$ with an operator $i: ring \rightarrow ring$ that satisfies $x \cdot i(x) \cdot x = x$. This function $i$ need not be unique, because $i(0)$ can take any value in $A$. However, if $j(x)$ is another function on the domain of $A$ such that
for all $x$, $x \cdot j(x) \cdot x = x$, then for all $x$, 

\begin{center}
$i(x) \cdot x \cdot i(x) = j(x) \cdot x \cdot j(x)$. 
\end{center}
To see this, write: $p(x) = i(x) \cdot x \cdot i(x)$ and $q(x) = j(x) \cdot x \cdot j(x)$. 
Now $x \cdot x \cdot p(x) = 
x \cdot x \cdot  i(x) \cdot x \cdot i(x)= x \cdot x \cdot  i(x) = x$ and $p(x) \cdot p(x) \cdot x = 
i(x) \cdot x \cdot i(x) \cdot i(x) \cdot x \cdot i(x) \cdot x = i(x) \cdot x \cdot i(x) \cdot i(x) \cdot x =
x \cdot i(x) \cdot i(x) = p(x).$ An application of Lemma \ref{uniqueinv} establishes that $p(x) = q(x)$ for all $x$. 
It follows that $p$ is independent of the choice of $i$. 

Then expand $A^{\prime}$ to the $\Sigma$ algebra
$A^{\prime\prime}$ by introducing an inverse operator as follows:
\begin{center}
$x^{-1} =p(x)= i(x) \cdot x \cdot i(x)$. 
\end{center}

We will show that both $\mathit{Ril}$ and $\mathit{Ref}$ are satisfied. For $\mathit{Ril}$ we make use of the 
equations just derived for $p(-)$ and find: $x \cdot x \cdot x^{-1} = x \cdot x \cdot p(x) = x$.

Now $\mathit{Ref}$ has to be established for the proposed inverse operator. In order to prove that
$(u^{-1})^{-1} = u$, write $x = u^{-1}$, $y = x^{-1}$ and $z = u$. 

Then, using straightforward calculations, we obtain: $x \cdot x \cdot y = x$, $y \cdot y \cdot x = y$, $x \cdot x \cdot z = x$ and $z \cdot z \cdot x = z$. It follows by Lemma \ref{uniqueinv}  that $y = z$, which is the required identity.

To see that the expansion is unique suppose that two unary functions $p(-)$ and $q(-)$ both satisfy $\mathit{Ref}$ and
$\mathit{Ril}$. Using Lemma \ref{SIP12} both functions satisfy $p(x \cdot y) = p(x) \cdot p(y)$ and 
$q(x \cdot y) = q(x) \cdot p(y)$, respectively. Given an arbitrary $x$ we find: 
$x \cdot x \cdot p(x) = x$ by assumption on $p(-)$. Applying $p(-)$ on both sides we find $p(x \cdot x \cdot p(x)) = p(x)$, which using  $\mathit{SIP2}$ implies $p(x) \cdot p(x) \cdot p(p(x)) = p(x)$. Then, using $\mathit{Ref}$ we have $p(x) \cdot p(x) \cdot x = p(x)$. Similarly we find $x \cdot x \cdot q(x) = x$ and 
 $q(x) \cdot q(x) \cdot x = q(x)$. By means of Lemma \ref{uniqueinv} this yields $p(x) = q(x)$.
\end{proof}

The uniqueness of inverse as an expansion of commutative rings satisfying $\mathit{Ref}$ and $\mathit{Ril}$ indicates that the inverse operation can be implicitly defined on a commutative von Neumann regular ring. The Beth definability theorem implies the existence of an explicit definition for inverse. In this case the application of Beth definability is inessential, however, because from the proof of Lemma \ref{vNrrexpansion} an explicit definition can be 
inferred for $y = x^{-1}$:
\[\exists z. (x \cdot z \cdot x = x\,  \&\,  y = z \cdot x \cdot z).\]

\section{The embedding theorem}
Because the theory of meadows is equational we know from universal algebra 
(see \cite{Meinke&Tucker.92, Wechler.92}) that:
\begin{theorem}\label{Birkhoff} The class of meadows is closed under subalgebras, direct products and homomorphic images. 
\end{theorem}
Thus, every subalgebra of a product of zero totalized fields is a meadow.
Our main task is to show that  
every non-trivial meadow is isomorphic to a subalgebra of a product of zero totalized
fields. First, we recall some basic properties of commutative rings, which can be found in many textbooks (e.g., \cite{McCoy.64}). 

\subsection{Preliminaries on rings}
 Let $R$ be a commutative ring. 
 An {\em ideal} in a ring $R$ is a subset $I$ with 0, and such that if
 $x,y\in I$ and $z\in R$, then $x+y\in I$, and $z\cdot x\in I$. $R$ itself
 and \{0\} are the trivial ideals. Any other ideal is a {\em proper
 ideal}.

 The ideal $R\cdot x=\{y\cdot x|~~y\in R~\}$ is  {\em the principal
 ideal} generated by $x$. Since $R$ has a unit, the generator
 $x=x\cdot 1$ is in $R\cdot x$. This is the smallest ideal that
 includes $x$.

If $I$ is an ideal then the following relation is a $\Sigma_{CR}$ congruence:
 \[x\equiv y\quad \mathit{iff}\quad x-y \in I.\]

The set of classes $R/I$ is a ring. The {\em quotient map} maps
every element $a$ of $R$ to its equivalence class, which is denoted
by $a+I$ or by $a/I$. The quotient map is a  $\Sigma_{CR}$ homomorphism from
$R$ onto $R/I$ (an epimorphism).
 It is clear what it means that $I$ is a maximal ideal in $R$.

\begin{lemma} Every ideal is contained in (at least one) maximal ideal.
\end{lemma}

\begin{proof} The union of a chain of ideals containing $I$ and
not 1 does not include 1. Therefore, by Zorn's lemma there is a maximal
such ideal.  
\end{proof}

\begin{lemma}  I is a maximal ideal iff $R/I$ is a field.
\end{lemma}

\begin{proof} If $x$ is not in $I$ then the ideal generated by $I$ and
$x$ is $R$. Hence for some $i$ in $I$ and $y$ in $R$ we have $1=i+xy$. It
follows that the classes of $x$ and of $y$ are inverse to each other.
Since $x$ is arbitrary outside $I$, every class except for the class 0
(i.e, the set $I$) has an inverse.
\end{proof}

Recall that 
$e\in R$ is called an {\em idempotent} if $e\cdot e=e$.
 \begin{propos}\label{idemp} Let $e\in R$ be an idempotent and $e\cdot R$
the principal ideal that it generates. Then\\
1. $e$ is a unit in the ring $ e\cdot R$,\\
2.  the mapping $H(a)=e\cdot a$ is a $\Sigma_{CR}$ homomorphism from $R$
onto the ring $e\cdot R$,\\
3.  For every $x\in R$:     $x\in e\cdot R ~~~\mathit{iff}~~~e\cdot x=x$.
\end{propos}

\begin{proof}~\\
1.  Note that $e=e\cdot 1$ and therefore $e\in e\cdot R$. For every element
$e\cdot a$ in $e\cdot R$ we have $e\cdot(e\cdot a)=e\cdot a$, by
associativity, and because $e\cdot e=e$. Therefore $e$ is a unit
in $e\cdot R$.\\
2. $H$ is a $\Sigma_{CR}$  homomorphism since:

 $e\cdot 0=0$ and $e\cdot 1=e$, so that zero is mapped to zero,
 and the unit is mapped to the unit.

 $e(a+b)=e\cdot a +e\cdot b$ and $e\cdot(-a)=-e\cdot a$, so that
 $+$ and $-$ are preserved.

 $e(f\cdot g)=(e\cdot e)
(f\cdot g)=(e\cdot f)(e\cdot g)$ so that multiplication is
preserved.\\
3. If $x\in e\cdot R$ then $e\cdot x=x$ by (1). And if $x=e\cdot
x$ then the right side testifies that it is an element of $e\cdot
R$.
\end{proof}

\subsection{Principal ideals in a meadow}
Let $R$ be a non-trivial meadow, and $x\in R$ a non zero element.
Note that by $Ril$, $1_x$ is an idempotent.

\begin{propos}
The principal ideal $x\cdot R$ has the following properties:\\
(a)  $1_x\cdot R=x\cdot R$, and  $x, 1_x$ and $x^{-1}$ are all in
$x\cdot R$.\\
(b) $x\cdot R$ is a ring with a unit, $x$ is invertible in the ring
and $H(y)=1_x\cdot y$ is a $\Sigma_{CR}$ homomorphism from $R$ onto $x\cdot
R$.
\end{propos}

\begin{proof}
(a)  Now $1_x=x^{-1}\cdot x$ hence  $1_x \in x \cdot R$, and $x=x\cdot 1_x$
hence $x \in 1_x \cdot R$. Therefore, $x\cdot R =1_x \cdot R$. Consequently, both
$x$ and $1_x$ belong to the ideal that they generate, and since
$x^{-1}=1_x\cdot x^{-1}$, $x^{-1}$ is also in $1_x\cdot R$.

(b) Since $1_x$ is an idempotent, this is Proposition \ref{idemp}. Note that $x$ is
invertible since $x\cdot x^{-1}$ is the unit in this ring, and $
x^{-1}$ is also in it.
\end{proof}

\begin{propos}\label{HOMr} Let $R$ be a meadow.
For every non-zero $x\in R$ there is a $\Sigma_{CR}$ homomorphism 
$H_x: R \rightarrow F_x$ from $R$ onto a
zero totalized field $F_x$ with $H_x(x)\neq0$.
\end{propos}

\begin{proof}
Let $x \neq 0$ be given, and let $I$ be a maximal ideal in the
ring $1_x \cdot R$. Then $R/I$ is a field, and the mapping
$H_x(y)=(y\cdot 1_x)/I$ is a $\Sigma_{CR}$ homomorphism as it is the composition of two
$\Sigma_{CR}$ homomorphisms. Now $H_x(x)=x/I$
and  $H_x(x)\neq 0$ because if an invertible element of $1_x \cdot R$ 
is mapped to $0$ by the quotient map, then $1=0$ in the
quotient $R/I$.
\end{proof}

\begin{propos}\label{HOMm}
If $H: R \rightarrow F$ is a $\Sigma_{CR}$ homomorphism from a meadow $R$ 
into a zero totalized field 
$F$ then $H$ preserves inverses and so is a $\Sigma$ homomorphism.
\end{propos}
\begin{proof}
If $H(x)=0$ then $H(1_x)=H(x\cdot x^{-1})=H(x)\cdot
H(x^{-1})=0$ so that also implies $H(x^{-1})=H(1_x\cdot
x^{-1})=H(1_x)\cdot H(x^{-1})= 0=H(x)^{-1}$.
 The latter holds because $F$ is 
zero totalized. Secondly, we consider the case
that $H(x)\neq 0$. Then $H(x)=H(1_x\cdot x)=H(1_x) \cdot H(x)$ which
proves that $H(1_x)=1$, by cancellation in fields. In other words
$1=H(x\cdot x^{-1})=H(x)\cdot H(x^{-1})$, which proves that
$H(x^{-1})=H(x)^{-1}$ using Proposition \ref{XYZ}.
\end{proof}

The image of $H$ is subfield of $F$, so it follows that given $R$ and non-zero $x \in R$ a 
meadow homomorphism onto a field $F$ can be found which maps $x$ to a 
non-zero element of $F$.
Using these preparations, we can prove the embedding theorem:
\begin{theorem}
 A $\Sigma$ structure is a non-trivial meadow if and only if it is a $\Sigma$-substructure
of a product of zero totalized fields.
\end{theorem}
\begin{proof}
By Theorem \ref{Birkhoff} a $\Sigma$ subalgebra of a  product of zero totalized fields
is always a meadow.

Let $R$ be a meadow. Combining Propositions \ref{HOMr} and \ref{HOMm}, for each nonzero $x$ in $R$ there is a field $F_x$ and a $\Sigma$ homomorphism
$H_x: R \rightarrow F_x$, such that $H_x(x)\neq 0$.

We define the product of fields:  $ K=\prod_{x\in R} F_x$. $K$
is a meadow with the operations defined at each coordinate. We define the
map $H$ from $R$ to the product as follows: for every $z$ in $R$, $H(z)$
is the vector that has $H_x(z)$ in the place $x$. Since $H_x$ is a
$\Sigma$-homomorphism with respect to all meadow operations, following the principles
of universal algebra, the same is true for $H$ as well.

If $z\neq 0$ then  $H_z(z)\neq 0$ and consequently $H(z)\neq 0$. Therefore $H$
is a $\Sigma$-monomorphism, which concludes the proof.
\end{proof}

\begin{corollary}
A finite non-trivial meadow $R$  is a $\Sigma$-substructure of a finite product of finite fields.
\end{corollary}

\subsection{Equational theory of  zero totalized fields}
The equational theory of zero totalized fields and of meadows are the
same. More precisely:
\begin{theorem}\label{Completeness}
For every $\Sigma$-equation $e$, 
$Alg(\Sigma, \mathit{ZTF}) \models e \Leftrightarrow Alg(\Sigma, \mathit{Md}) \models e.$
\end{theorem} 

\begin{proof}
Let $e$ be an equation that holds  in every zero totalized field, then it 
holds also in every product of fields and in every $\Sigma$ subalgebra of
a product of fields, and therefore, by the embedding theorem, also
in every non-trivial meadow. Evidently, every equation holds in the trivial meadow as well.

The other way around, that equations true for all meadows hold in all zero totalized 
fields, is obvious because zero totalized fields are a subclass of meadows.
\end{proof}

\subsection{Conditional equational theory of  zero totalized fields}
As an application of Theorem \ref{Completeness}, we prove a stronger result, namely: the conditional equational theories of zero totalized fields and of meadows are the
same. More precisely:
\begin{theorem}\label{CCompleteness}
For every conditional $\Sigma$-equation $e$, 
$Alg(\Sigma,$
\[ \mathit{ZTF}) \models e \Leftrightarrow Alg(\Sigma, \mathit{Md}) \models e.\]
\end{theorem} 

\begin{proof}
Let $t^1_1 = t^1_2\, \&\, \ldots \,\&\, t^i_1 = t^i_2\, \& \,\ldots \,\& \,t^n_1 = t^n_2 \rightarrow  t_1 = t_2$ be a conditional
equation that holds  in every zero totalized field. Without loss of generality, it may be assumed that each right-hand side equals 0, using $r = s \Leftrightarrow r-s = 0$. So we assume that 
\\
$t_1 = 0\, \&\, \ldots \,\&\, t_i=0\, \& \,\ldots \,\& \,t_n= 0 \rightarrow  t = 0$ holds in all zero totalized fields. If $n=0$ the case reduces to that of equations and the conclusion follows from Theorem \ref{Completeness}. Let the $\Sigma$ term $C(-,-)$ be given by 
\[C(x,y) = (1-\frac{x}{x}) \cdot y.\] Now, by inspection of zero totalized fields, one has:
\[Alg(\Sigma, \mathit{ZTF}) \models t_1 = 0 \rightarrow t = 0 
\Leftrightarrow Alg(\Sigma, \mathit{ZTF}) \models C(t_1,t)=0.\]
As a consequence, $Alg(\Sigma, \mathit{Md}) \models C(t_1,t)=0$. Now, 
$\mathit{Md} \cup \{C(t_1,t)=0\} \vdash t_1 = 0 \rightarrow t = 0$ and consequently 
$\mathit{Md} \vdash t_1 = 0 \rightarrow t = 0$ and, of course, $\mathit{Md} \models t_1 = 0 \rightarrow t = 0$.

In the case of $n=2$ we assume that all zero totalized fields satisfy $t_1 = 0 \, \& \, t_2 = 0 \rightarrow t = 0$. We will make use of the following fact which holds in all meadows:
\[x = 0\, \&\, y = 0 \Leftrightarrow \frac{x \cdot y}{x \cdot y} - \frac{x}{x} - \frac{y}{y} = 0\]
Here ``$\Rightarrow$'' is immediate and to see ``$\Leftarrow$'' multiply both sides with $ x$ thus obtaining:
\[ \frac{ x \cdot x \cdot y}{x \cdot y} - \frac{x \cdot x}{x} - \frac{ x \cdot y}{y} = x \cdot 0\]
and, using $\mathit{Md}$,
\[ \frac{x  \cdot y}{ y} - x  - \frac{x \cdot y}{y} = 0\]
which implies $x=0$. Similarly, one derives $y = 0$. We write 
$U(x,y) = \frac{x \cdot y}{x \cdot y} - \frac{x}{x} - \frac{y}{y}$. Now using $U(x,y) = 0 \Leftrightarrow x=0 \, \& \, y = 0$, we find:
\[Alg(\Sigma, \mathit{ZTF}) \models t_1 = 0 \, \& \, t_2 = 0 \rightarrow t = 0 
\Leftrightarrow Alg(\Sigma, \mathit{ZTF}) \models C(U(t_1,t_2),t)=0.\]
Using Theorem \ref{Completeness}, we find that $\mathit{Md} \models C(U(t_1,t_2),t)=0$ and, from this fact using the known properties of $U(-)$ and $C(-,-)$, one easily derives 
$\mathit{Md} \models t_1 = 0 \, \& \, t_2 = 0 \rightarrow t = 0$. The cases $n =3, \dots$ require a repeated nested use of $U(-)$. The straightforward details have been omitted and we only illustrate the encoding of conditional equations into equations in the case $n=3$:
\[Alg(\Sigma, \mathit{ZTF}) \models (\bigwedge_{i=1}^{i=3}t_i = 0) \rightarrow t = 0 
\Leftrightarrow Alg(\Sigma, \mathit{ZTF}) \models C(U(U(t_1,t_2),t_3),t)=0.\]
\end{proof}

\section{Finite meadows}
 
 As usual, we will define $\underline{0}$ as $0$ and $\underline{k+1}= \underline{k}+1$.
The characteristic of a meadow is the smallest natural number $k \in \nat$ such that $k > 0$ and $\underline{k}=0$.
The equation $\underline{k}=0$ will be referred to as $Z_{k}$.
We recall that a natural number $k$ is called {\em squarefree} if its prime factor decomposition is the product of distinct primes.
\begin{lemma}\label{DistinctPrimes}
Let $M$ be a meadow of finite characteristic $k >0$. Then $k$ is squarefree.
\end{lemma}
\begin{proof}
Let $M \models \underline{k} = 0$. Suppose $k$ has two repeated prime factors, $k= p \cdot p \cdot q$. Then, using $Ril$ we have
\begin{center}
$\underline{p}\cdot \underline{q} = (\underline{p} \cdot \underline{p} \cdot \underline{p}^{-1}) \cdot \underline{q} = (\underline{p} \cdot \underline{p} \cdot \underline{q}) \cdot \underline{p}^{-1} =  \underline{k} \cdot \underline{p}^{-1} =0 \cdot \underline{p}^{-1} = 0$.
\end{center}
Thus, $k$ is not the characteristic which is a contradition.
\end{proof}

Thus, from Lemma \ref{DistinctPrimes}, the possible finite characteristics have the form $k=p_1 \ldots p_n$ where the $p_i$ are all distinct primes. All finite meadows have finite characteristic. It follows that if a finite meadow $M$ consists of an initial segment of the numerals \underline{0}, \ldots, \underline{k-1} (like the prime fields of positive characteristic) its cardinality $\#(M)=k$ can only be a product of different primes.

\begin{defin}
Let $\mathit{Md}_{k}$  be the initial algebra of $\mathit{Md} \cup \{ Z_{k}\}$.  
\end{defin}

What are the initial algebras? Clearly, $\mathit{Md}_{k}$ has finite characteristic $\leq k$. Notice the following:

\begin{lemma}
If $l$ divides $k$ then the $\mathit{Md} + Z_{l} \vdash  Z_{k}$. Thus, if $l$ divides $k$ then there is a $\Sigma$ epimorphism $ \phi \colon  \mathit{Md}_{k} \to \mathit{Md}_{l}$, i.e.,  $\mathit{Md}_{l}$ is a homomorphic image of $\mathit{Md}_{k}$.
\end{lemma}
 
Thus, we have that for $k=p_1 \ldots p_n$ where the $p_i$ are all distinct primes we have a
$\Sigma$ epimorphism $\phi \colon  \mathit{Md}_{k} \to \mathit{Md}_{p_i}$. Furthermore, it can be seen that for $p$ a prime number,  $\mathit{Md}_{p}$ is the zero totalized prime field $\integers_p$ of characteristic $p$. To see this notice that for each $x$ different from $0$ there is an $y$ with $x \cdot y = 1$. It follows that the zero totalized prime field mod $p$ satisfied \textit{Iel} (see Section~\ref{Iel})
and for that reason it is a meadow. As a consequence we have a $\Sigma$ epimorphism 
$ \phi \colon  \mathit{Md}_{k} \to \integers_{p_i}$.

\begin{theorem}\label{DistinctPrimes2}
If $k$ is squarefree then $\mathit{Md}_{k}$ has $k$ elements.
\end{theorem}
\begin{proof}
If $k=p_1 \ldots p_n$ is a product of different primes that is no prime factor appears twice then we first show that $\mathit{Md}_{k}$ has at least $k$ elements. To see this notice that for each prime factor $p$ of $k$ the prime field $\integers_p$ of characteristic $p$ is a model of $\mathit{Md}_{k}$ (as the equation $Z_{p}$ implies $Z_{k}$). Because that structure is a quotient of the additive group of $\mathit{Md}_{k}$ its number of elements is a divisor of the cardinality $\#(\mathit{Md}_{k})$ of $\mathit{Md}_{k}$. As a consequence $\#(\mathit{Md}_{k})$ is a multiple of all factors of $k$ and because $k$ contains all of them only once $\#(\mathit{Md}_{k}) \geq k$.

In order to prove that $\#(\mathit{Md}_{k}) = k$ it suffices to find an inverse (in the sense of a meadow) for each $\underline{n}$ for $n<k$ of the form $\underline{m}$ for $m < k$. We may assume that $k>0$ otherwise the inverse is obvious. To find the inverse consider the power series $\underline{n}^{0} (=1), \underline{n}^{1} , \underline{n}^{2} ..$. Each value in this series is of the form $\underline{m}$ for $m < k$ because arithmetic is done modulo $k$. Therefore there are $k$ and $l$ with $k>l+1>0$ such that $\mathit{Md}_{k}\models  \underline{n}^{k} = \underline{n}^{l}$. Let $k-1-l=i$. Notice that $i \geq 0$. Working in $\mathit{Md}_{k}$ by $\mathit{SIP2}$
 we have $\underline{n}^{-k} = \underline{n}^{-l}$, and thus $\underline{n}^{-1} = \underline{n}^{-k} \cdot \underline{n}^{k-1}= \underline{n}^{-l}\cdot \underline{n}^{k-1}= \underline{n}^{k-1-l}=\underline{n}^{i}$. This demonstrates that the inverse is a numeral (modulo $k$) as required.
\end{proof}

It follows from the proof that the interpretation of inverse is unique in a minimal finite meadow.
Recall that an algebra is minimal when it has no subalgebras or, equivalently,  is generated by elements named in its signature. By Lemma \ref{DistinctPrimes2}, if $k$ is a product of different primes then $\mathit{Md}_{k}$ is the minimal meadow of characteristic $k$. 
It also follows from the proof that $\mathit{Md}_{k}$ consists of $0,\ldots k-1$.\\

\noindent {\bf Example 1.} Concrete examples can be easily given, for instance $\mathit{Md}_6$ has the following inverse function: $0^{-1} = 0, 1^{-1} = 1, 2^{-1} = 2, 3^{-1} = 3, 4^{-1} = 4,$ and $ 5^{-1} = 5$. $\mathit{Md}_6$ is the smallest non-trivial minimal meadow which is not a field.\\ 

\noindent {\bf Example 2.} In $\mathit{Md}_{10}$ the inverse function is given by: $0^{-1} = 0, 1^{-1} = 1, 2^{-1} = 8, 3^{-1} = 7, 4^{-1} = 4, 5^{-1} = 5, 6^{-1} = 6, 7^{-1} = 3, 8^{-1} = 2,$ and $ 9^{-1} = 9$.\\

\noindent {\bf Example 3.} Consider $\mathit{Md}_{4}$. This is a non-minimal meadow because its size of four elements exceeds its characteristic.  The inverse function is the identity function. 
$\mathit{Md}_{4}$ is the smallest non-trivial meadow which is not a field.\\

\begin{lemma}\label{monomorphism}
Let $M$ be a meadow of finite characteristic $k >0$. Then there is a $\Sigma$ monomorphism  $ \psi \colon  \mathit{Md}_{k} \to M$.
\end{lemma}
\begin{proof}
If $M$ has characteristic $k$ then $M \models \underline{k}=0$. Thus, by initiality, there is a 
$\Sigma$ homomorphism $ \psi \colon  \mathit{Md}_{k} \to M$. If this map were not injective then $M$ would have characteristic lower then $k$.
\end{proof}

\begin{lemma}\label{isomorphism}
Let $M$ be a minimal meadow of finite characteristic $k >0$. Then  $\mathit{Md}_{k}$ and $M$
are $\Sigma$ isomorphic.
\end{lemma}
\begin{proof}
If $M$ has characteristic $k$ then $M \models \underline{k}=0$. Thus, following the previous lemma there is a $\Sigma$ monomorphism $ \psi \colon  \mathit{Md}_{k} \to M$. Because $M$ is minimal, $\psi$ is surjective as well.
\end{proof}

\begin{lemma}\label{xx}
Let $M$ be a meadow of prime cardinality $p$. Then  $M$ is the zero totalized prime field of cardinality $p$.
\end{lemma}
\begin{proof}
If $M$ has characteristic $k$ then $k > 0$ is the cardinality of the smallest additive subgroup of
$M$ which contains $1$. Thus $k$ divides $p$ and hence $k = p$ which implies that $M$ is minimal. Following Lemma~\ref{isomorphism}
$\mathit{Md}_{k}$ is isomorphic with $M$. At the same time the zero totalized prime field of cardinality $p$ is a meadow and according to 
Lemma \ref{isomorphism}  it is also isomorphic to $\mathit{Md}_{k}$. \end{proof}

\begin{lemma}\label{survey}
All finite and minimal meadows are of the form $\mathit{Md}_{k}$ for some positive natural number
 $k$. 
\end{lemma}

\begin{proof}
Let $M$ be a finite meadow. Then $M$ has a finite characteristic, say $k$. 
By Lemma~\ref{isomorphism}, there is an isomorphism $ \psi \colon  \mathit{Md}_{k} \to M.$
\end{proof}

If its non-zero characteristic is not a prime, a finite meadow has proper zero-divisors and fails to be an integral domain and, of course, it is no field either.

\begin{lemma}
If $k = p_{1}^{\alpha_{1}}  \ldots  p_n^{\alpha_{n}}$ then $\mathit{Md}_{k} \cong \mathit{Md}_{p_1 \ldots p_n}$. Therefore, if $k$ and $l$ have the same set of prime factors then $\mathit{Md}_{k} \cong \mathit{Md}_{l}$. 
\end{lemma}

\begin{proof} Using the same argument as in Lemma \ref{DistinctPrimes}, we can show that for 
 $p_1, \ldots,  p_n$ any primes and $k = p_{1}^{\alpha_{1}}  \ldots  p_n^{\alpha_{n}}$ we have $\mathit{Md}_{k} \cong \mathit{Md}_{p_1 \ldots p_n}$.
Suppose that 
$k = p_{1}^{\alpha_{1}}  \ldots  p_n^{\alpha_{n}} $ and $l = p_{1}^{\beta_{1}}  \ldots  p_n^{\beta_{n}} $. Then by the first part of the lemma, 
$\mathit{Md}_{k} \cong \mathit{Md}_{p_1 \ldots p_n}$ and $\mathit{Md}_{l} \cong \mathit{Md}_{p_1 \ldots p_n}$
and hence $\mathit{Md}_{k} \cong \mathit{Md}_{l}$.
 
\end{proof}

\section{Concluding remarks and further questions}
We notice that a conference version of this paper, though with a quite different emphasis of presentation,  has appeared as \cite{BHT07}.

The theory of meadows depends upon the {\it formal} idea of a total inverse operator. We do not claim that division by zero is possible in numerical calculations involving the rationals or reals.  But we do 
claim that zero totalized division is logically, algebraically and computationally useful: for some applications, allowing zero totalized division in formal calculations, based on equations and rewriting, is appropriate because it is conceptually and technically simpler than the conventional concept of partial division. Furthermore, one can make arrangements to track the use of the inverse operation in formal calculations and classify them them as safe or unsafe dependent upon $0^{-1}$ is invoked: see \cite{Bergstra&Tucker.07}. We expect these areas to include elementary school algebra, specifying and understanding gadgets containing calculators, spreadsheets, and declarative programming. Of course, further research is necessary to test these expectations: at present, our theory of meadows is a theory of zero totalized division, constitutes a generalization of the theory of fields, and is known to be useful in specifying numerical data types using equations.

There are many opportunities for the further development of the theory of meadows: logically, algebraically, and through applications. Consider some computational and logical open questions that add to the questions posed in \cite{Bergstra&Tucker.05a}:

Is the equational theory of meadows decidable? Is its conditional equational theory decidable?

Does $\mathit{Md}$, or a useful extension of it, admit Knuth-Bendix completion?

Returning to the equational theory of meadows, following \cite{Bergstra&Tucker.05a}, let
 $Z(x) = 1- x \cdot x^{-1}$. For $n >0$, let $L_n$ be the equation: $Z(1 +x_1^2 +....+x_n^2) = 0$. Clearly from $CR$ it follows that $L_k$ implies $L_n$ when $k>n$. 
All $L_n$ are valid in the zero totalized field of rational numbers. From \cite{Bergstra&Tucker.05a} 
and Proposition \ref{SIP12}, it follows that $\mathit{Md} + L_4$ constitutes an initial algebra specification of the zero totalized field of rational numbers, which indicates the relevance of $L_4$. Now, conversely, the question arises if $\mathit{Md} + L_n$ proves $L_k$ (again assuming $k>n$). 

A related problem is to characterize the initial algebras of $\mathit{Md} + L_n$  for $n=1$, 
$n=2$, and $n=3$. It is easy to see that $\mathit{Md} + L_1$ is not a specification of the rationals because it is satisfied by the prime field of characteristic three, which is not a homomorphic image of the 
initial algebra of $\mathit{Md} + L_4$. 

A restricted version of Theorem \ref{Completeness} for equations between closed terms 
only, was 
shown in \cite{Bergstra&Tucker.05a}. That proof is longer and more syntactic in style and uses 
 a normal form result and straightforward induction, in spite of the fact that the result is weaker. 
However, it provides the additional information that the initial algebra of $\mathit{Md}$ is a 
computable algebra. 
The proof given here uses the maximal ideal theorem, which is weaker than the 
axiom of choice, but still independent of the axiom system $\mathit{ZF}$ for set theory. The use of maximal ideals provides a  simple and readable proof. In \cite{BP08}, however, a proof is given in the proof theoretic style. That proof is more general and it provides the information that the equational consequences of $\mathit{Md} + L_4$ coincide with the equations valid in all zero-totalized fields that satisfy $L_4$, which seems not to follow from a proof using maximal ideals.
 
Finally, let us note that questions may emerge from the perspective of pure algebra, where the properties of invertibility and symmetry are central.  The representation results here are closely related to early results on subdirect products of rings of McCoy \cite{McCoy.38} and Birkhoff \cite{Birkhoff.44}. 

The results leading up to the representation and completeness theorems may be investigated for non-commutative rings. The theory of von Neumann regular rings is primarily about non-commutative rings.  As is always the case, the transition from commutative to non-commuutative rings is a delicate operation, leading to a ramification of properties. In \cite{Bergstra&Hirshfeld&Tucker.07} we have isolated a number of concepts and proved generalizations of the main results here to skew fields and skew meadows.

We define a {\em skew meadow} to be an expansion of a non-commutative ring with an inverse operator that satisfies these two equations:
\begin{eqnarray}
(x^{-1})^{-1}& = & x \label{Ref} \\
x \cdot (x \cdot x^{-1}) & = & x \label{Ril}
\end{eqnarray}
Thus, the equations for skew meadows result from the equations for meadows, by simply dropping commutativity of multiplication and including a second distributivity law: a {\em meadow} is a commutative skew meadow.  Actually, the simplicity of this generalisation is a technical achievement for there are several interesting equations that are equivalent in the commutative case but in differ in the non-commutative case; also, these equations must be distinguished as rewrite rules.  In \cite{Bergstra&Hirshfeld&Tucker.07} we consider several related types of non-commutative ring.

\end{document}